\documentclass[graybox]{svmult}

\usepackage{type1cm}        
\usepackage{makeidx}         
\usepackage{graphicx}        
\usepackage{multicol}        
\usepackage[bottom]{footmisc}

\usepackage{amsmath}
\renewcommand\vec[1]{%
    \mathbf{\itshape #1}
}

\newcommand\trasp[1]{%
    {#1}^\mathsf{T}
}

\DeclareMathOperator{\de}{d\!}
\renewcommand{\vec}[1]{\bm{\mathrm{#1}}}
\newcommand{\up}[1]{\ensuremath{\mathrm{#1}}}
\newcommand{\od}[2]{\frac{\de{#1}}{\vphantom{l^l}\de{#2}}}
\newcommand{\abs}[1]{\left|#1\right|}

\newcommand{\Qe}{\vec{V}_\up{e}}
\newcommand{\Qs}{{\vec{V}^\ast}}
\newcommand{\Vs}{{\vec{V}^\ast}}
\newcommand{\Js}{{\vec{J}^\ast}}
\newcommand{\Ss}{{\vec{S}^\ast}}
\newcommand{\Bs}{{\vec{B}^\ast}}
\newcommand{\Cs}{{\vec{C}^\ast}}
\newcommand{\tss}{{t^\ast}}

\makeindex             

\title*{A solver for stiff finite-rate relaxation in Baer--Nunziato two-phase flow models}
\titlerunning{A solver for stiff finite-rate relaxation in Baer--Nunziato two-phase models}
\author{Simone Chiocchetti \and Christoph M\"uller}
\institute{
Simone Chiocchetti \at Laboratory of Applied Mathematics, University of Trento, via Mesiano 77, 38123 Trento, Italy, \email{simone.chiocchetti@unitn.it},\\
Christoph M\"uller \at Institute of Aerodynamics and Gasdynamics, Pfaffenwaldrig 21, 70569 Stuttgart, Germany.}

\begin{document}

\maketitle
\abstract{ 
In this paper we present a technique for constructing robust solvers for stiff algebraic source terms, 
such as those typically used for modelling relaxation processes in 
hyperbolic systems of partial differential equations describing two-phase flows, namely models of the
Baer--Nunziato family. The method is based on an exponential integrator 
which employs an approximate linearised source term operator 
that is constructed in such a way that one can
compute solutions to the linearised equations avoiding 
any delicate matrix inversion operations.}

\section{Introduction}
Stiff algebraic source terms, accounting for mechanical relaxation and phase transition 
in two-phase flow models of the Baer--Nunziato type \cite{baernunziato, pelanti, spb}, 
are one of the key difficulties in computing solutions to these systems of hyperbolic partial differential 
equations (PDE). Their accurate solution is relevant for the study of droplet 
dynamics with Baer--Nunziato models. These weakly compressible phenomena
can be accurately described by the reduced models that assume instantaneous pressure
and velocity equilibrium like the one forwarded by Kapila et al. \cite{kapila}. Solving 
more general sets of equations like \cite{baernunziato, pelanti, spb} in the stiff relaxation limit
gives results that are similar to those obtained from the instantaneous equilibrium model, while allowing
more modelling flexibility, since less physical assumptions have to be made.

A simple computational strategy for dealing with stiff sources is the \emph{splitting} 
approach \cite{strang,fractionalstep}.
The procedure consists of two steps: at each timestep, first one solves the homogeneous part of the PDE 
\begin{equation} \label{eq:modelpde}
    \partial_t \vec{Q} + \nabla\cdot\vec{F}(\vec{Q}) + \vec{B}(\vec{Q})\,\nabla\vec{Q} = \vec{S}(\vec{Q}),
\end{equation}
for example with a 
path-conservative \cite{castro2006,pares2006} MUSCL--Hancock \cite{muscl} method, 
obtaining a preliminary solution $\vec{Q}_\up{H}$ and then
one can use this state vector as initial condition for the Cauchy problem
\begin{equation} \label{eq:modelcauchy}
\od{\vec{Q}}{t} = \vec{S}(\vec{Q}),\quad \vec{Q}(t_n) = \vec{Q}_\up{H},\quad t\in ({t_n,\  t_{n+1}}),\\
\end{equation}
of which the solution will then yield the updated quantities at the new time level $t_{n+1}$.
This way, the problem is reduced to the integration of a system of ordinary differential equations (ODE), 
and general-purpose ODE solvers or
more specialised tools can be employed for this task.

It is often the case that the time scales associated with relaxations sources 
are much shorter than those given by the stability condition of the PDE scheme, 
thus one must be able to deal with source terms that are potentially stiff.
In order to integrate stiff ODEs with conventional explicit solvers,
one has to impose 
very severe restrictions on the maximum timestep size,
and for this reason implicit methods are commonly
preferred \cite{torobook}. Unfortunately implicit solvers, are, on a per-timestep basis, much more
expensive than explicit integrators, and they still might require variable sub-timestepping in order to
avoid under-resolving complex transients in the solution. 

In this work, we will develop a technique for constructing a solver 
for stiff finite-rate mechanical 
relaxation sources, specifically those encountered in models of the Baer--Nunziato type.

The proposed method overcomes the issues typical of explicit solvers with three concurrent strategies:
first, the update formula is based on exponential integration \cite{expint1, expint2}, in order to mimic at the discrete level
the behaviour of the differential equation; second, 
information at the new 
time level $t_{n+1}$ is taken into account by iteratively updating a linearisation
of the ODE system, this is achieved without resorting to a fully implicit method like those introduced in \cite{butcher1964}, 
and for which
one would need to solve a system of nonlinear algebraic 
equations at each timestep $t_n$; third, the method incorporates a simple and effective adaptive timestepping
criterion, which is crucial for capturing abrupt changes in the state variables and dealing with
the different time scales that characterise the equations under investigation.

\section{Model equations}
We are interested in the solution of two-phase flow models of the Baer--Nunziato family, 
which can be written in the
general form \eqref{eq:modelpde},
with a vector of conserved variables defined as
\begin{equation}
    \vec{Q} = \trasp{[\alpha_1\,\rho_1,\ \alpha_2\,\rho_2,\  \alpha_1\,\rho_1\,\vec{u}_1,\ \alpha_2\,\rho_2\,\vec{u}_2,\ 
\alpha_1\,\rho_1\,E_1,\ \alpha_2\,\rho_2\,E_2,\  \alpha_1]},
\end{equation}
a conservative flux $\vec{F}$ and a non-conservative term $\vec{B}\,\nabla\vec{Q}$ written as
\begin{equation}
    \vec{F}(\vec{Q}) = \begin{bmatrix}
    \alpha_1\,\rho_1\,\vec{u}_1\\
    \alpha_2\,\rho_2\,\vec{u}_2\\
    \alpha_1\left(\rho_1\,\vec{u}_1\otimes\vec{u}_1 + p_1\,\vec{I}\right)\\
    \alpha_2\left(\rho_2\,\vec{u}_2\otimes\vec{u}_2 + p_2\,\vec{I}\right)\\
    \alpha_1\,\left(\rho_1\,E_1 + p_1\right)\,\vec{u}_1\\
    \alpha_2\,\left(\rho_2\,E_2 + p_2\right)\,\vec{u}_2\\
    0
    \end{bmatrix},\quad
    \vec{B}(\vec{Q})\,\nabla\vec{Q} = \begin{bmatrix}
    0\\
    0\\
    - p_\up{I}\,\nabla\alpha_1 \\
    + p_\up{I}\,\nabla\alpha_1 \\
    - p_\up{I}\,\vec{u}_\up{I}\cdot\nabla\alpha_1\\
    + p_\up{I}\,\vec{u}_\up{I}\cdot\nabla\alpha_1\\
    \vec{u}_\up{I}\cdot\nabla\alpha_1
    \end{bmatrix},
\end{equation}
and a source term vector written as
\begin{equation} \label{eq:source}
    \vec{S}(\vec{Q}) = \begin{bmatrix}
        0\\
        0\\
        \lambda\left(\vec{u}_2 - \vec{u}_1\right)\\
        \lambda\left(\vec{u}_1 - \vec{u}_2\right)\\
        \lambda\left(\vec{u}_2 - \vec{u}_1\right)\cdot\vec{u}_\up{I} + \nu\,p_\up{I}\left(p_2 - p_1\right)\\
        \lambda\left(\vec{u}_1 - \vec{u}_2\right)\cdot\vec{u}_\up{I} + \nu\,p_\up{I}\left(p_1 - p_2\right)\\
        \nu\left(p_1 - p_2\right)
    \end{bmatrix}.
\end{equation}
Here we indicate with $\alpha_1$ and $\alpha_2$ the volume fractions of the first phase and of the second phase respectively, 
with $\rho_1$ and $\rho_2$ the phase densities, 
$\vec{u}_1 = \trasp{[u_1,\ v_1,\ w_1]}$ and 
$\vec{u}_2 = \trasp{[u_2,\ v_2,\ w_2]}$ indicate the velocity vectors,
$\alpha_1\,\rho_1\,E_1$ and $\alpha_2\,\rho_2\,E_2$ are the partial energy densities. The pressure fields are
denoted with $p_1$ and $p_2$, and the interface pressure and velocity are 
named $p_\up{I}$ and $\vec{u}_\up{I} = \trasp{[u_\up{I},\ v_\up{I},\ w_\up{I}]}$.
Finally, the parameters $\lambda$ and $\nu$ control the time scales for friction and pressure relaxation kinetics respectively.

In the following, we will study the system of ordinary differential equations arising from the source term~\eqref{eq:source} only,
that is, the one constructed as given in equation \eqref{eq:modelcauchy}
and specifically its one-dimensional simplification in terms of the primitive 
variables $\vec{V} = \trasp{\left[u_1,\ u_2,\ p_1,\ p_2,\ \alpha_1\right]}$, 
with an initial
condition $\vec{V}_0 = \trasp{\left[u_1^0,\ u_2^0,\ p_1^0,\ p_2^0,\ \alpha_1^0\right]}$. Since no source is
present in the mass conservation equations, they have a trivial solution, that is, $\alpha_1\,\rho_1$ 
and $\alpha_2\,\rho_2$ remain constant in time; for compactness, these quantities will be included in our 
analysis as constant parameters, rather than as variables of the ODE system.

The one-dimensional ODE system is written as
\begin{align}
            & \od{u_1}{t} = {\frac{\lambda}{\alpha_1\,\rho_1}}\,(u_2 - u_1)\label{eq:bn7odeu1}, \\[2mm]
            & \od{u_2}{t} = {\frac{\lambda}{\alpha_2\,\rho_2}}\,(u_1 - u_2)\label{eq:bn7odeu2}, \\[2mm]
            & \od{p_1}{t} = {\frac{\nu\,(p_I + k_{1a}\,p_1 + k_{1b})}{\alpha_1\,k_{1a}}}\,(p_2 - p_1) + 
                {\frac{\lambda\,(u_I - u_1)}{\alpha_1\,k_{1a}}}\,(u_2 - u_1)\label{eq:bn7odep1}, \\[2mm]
            & \od{p_2}{t} = {\frac{\nu\,(p_I + k_{2a}\,p_2 + k_{2b})}{\alpha_2\,k_{2a}}}\,(p_1 - p_2) + 
                {\frac{\lambda\,(u_I - u_2)}{\alpha_2\,k_{2a}}}\,(u_1 - u_2)\label{eq:bn7odep2},\\[2mm]
            & \od{\alpha_1}{t} = {\nu}\,(p_1 - p_2)\label{eq:bn7odea1}.
\end{align}

The choices for interface pressure and velocity are $p_I = p_2$ and $u_I = u_1$.
Finally, one can verify that, using the stiffened-gas equation of state for both phases, we have
$k_{1a} = 1/(\gamma_1 - 1)$, 
$k_{2a} = 1/(\gamma_2 - 1)$, 
$k_{1b} = \gamma_1\,\Pi_1/(\gamma_1 - 1)$, and
$k_{2b} = \gamma_2\,\Pi_2/(\gamma_2 - 1)$.
\section{Description of the numerical method}
\label{sec:method}

The methodology is described in the following with 
reference to a generic nonlinear first order Cauchy problem

\begin{equation} \label{eq:exactcauchy}
\od{\vec{V}}{t} = \vec{S}(\vec{V},\ t),\qquad
\vec{V}(t_n) = \vec{V}_n,
\end{equation}
for which the ODE can be linearised about a given state $\Qs$ and time $\tss$ as
\begin{equation}
\od{\vec{V}}{t} = \Bs + \Js(\Qs,\ \tss)\,(\vec{V} - \Qs).
\end{equation}
Here we defined the Jacobian matrix of the source ${\Js = \vec{J}(\Qs,\ \tss)}$ and analogously the source
vector evaluated at the linearisation state is ${\Bs = \vec{S}(\Qs,\ \tss)}$. 
We then introduce the vector 
\begin{equation}
\Cs = \Cs(\Bs,\ \Js) = \Cs(\Qs,\ \tss),
\end{equation}
which will be used as an indicator for the adaptive timestepping algorithm and may be
constructed for example listing all of the components of the matrix $\Js$ 
together with all the components of the vector $\Bs$ and the state $\Qs$, or only with a selection of
these variables, or any other relevant combination of the listed variables, that is, any group 
indicative of changes in the nature or the magnitude
of the linearised source operator.

It is then necessary to compute an accurate
analytical solution of the non-homogeneous linear Cauchy problem
\begin{equation} \label{eq:linearcauchy}
\od{\vec{V}}{t} = \Ss(\vec{V};\ \Qs,\ \tss) = \Bs + \Js(\Qs,\ \tss)\,(\vec{V} - \Qs), \quad\
 \vec{V}(t_n) = \vec{V}_n.\\
\end{equation}
We will denote the analytical solution of the IVP~\eqref{eq:linearcauchy} as $\Qe(t;\ \Ss,\ t_n,\ \vec{V}_n)$. 
As for $\Ss(\vec{V};\ \Qs,\ \tss)$, the semicolon separates the variable on which $\Qe$ and $\Ss$ continuously
depend ($t$ or $\vec{V}$) from the parameters used in the construction of the operators. 
The state vector at a generic time level $t_n$ is written as $\vec{V}_n$, the variable timestep size
is $\Delta t^n = t_{n+1} - t_n$.
\subsection{Timestepping}
Marching from a start time $t_0$ to an end time $t^\up{end}$ is carried out as follows.
First, an initial timestep size $\Delta t^0$ is chosen,
then, at each time iteration, the state $\vec{V}_{n+1}$ at the new time level $t_{n+1}$ is computed by means of the iterative 
procedure described below. The iterative procedure will terminate by computing a value for
$\vec{V}_{n+1}$, together with a new timestep size ${\Delta t^{n+1} = t_{n+2} - t_{n+1}}$ based on an estimator
which is embedded in the iterative solution algorithm. There is also the possibility that, 
due to the timestep size $\Delta t$ being too large, the value of $\vec{V}_{n+1}$ be flagged as not acceptable.
In this case, the procedure will return a new shorter timestep size for the current timestep 
${\Delta t^{n} = t_{n+1} - t_{n}}$ and a new attempt at the solution for $\vec{V}_{n+1}$ will
be carried out. Specifically, in practice we choose the new timestep size
to be half of the one used in the previous attempt.
\subsection{Iterative computation of the timestep solution}
At each iteration (denoted by the superscript $k$) we define an average state vector 
$\Qs_{n+1/2}^k = (\vec{V}_n + \Qs_{n+1}^{k-1})/2$
to be formally associated with an intermediate time level
$t_{n+1/2} = \left(t_n + t_{n+1}\right)/2$.
For the first iteration we need a guess value for $\Qs_{n+1}^{k-1}$, with the simplest choice being
$\Qs_{n+1}^{k-1} = \vec{V}_n$. Then the coefficients $\Cs_{n+1/2}^k$ are computed as
\begin{equation}
\Cs_{n+1/2}^k = \Cs_{n+1/2}^k(\Qs_{n+1/2}^k,\ t_{n+1/2}).
\end{equation}
In a joint way, one can build the affine source operator
\begin{equation}
\Ss_{n+1/2}^k = \Ss_{n+1/2}^k(\vec{V};\ \Qs_{n+1/2}^k,\ t_{n+1/2}).
\end{equation}
Then one can solve analytically
\begin{equation} \label{eq:approxcauchy}
\od{\vec{V}}{t} = \Ss_{n+1/2}^k(\vec{V};\ \Qs_{n+1/2}^k,\ t_{n+1/2}),\qquad
 \vec{V}(t_n) = \vec{V}_n,
\end{equation} 
by computing 
\begin{equation}
\Qs_{n+1}^k = \vec{V}_\up{e}\left(t_{n+1};\ \Ss_{n+1/2}^k,\ t_n,\ \vec{V}_n\right).
\end{equation}
It is then checked that the state vector $\Qs_{n+1}^k$ be physically admissible: in our case this means
verifying that internal energy of each phase be positive and that the volume fraction be bounded between 0 and 1.
Also one can check for absence of floating-point exceptions.
Additionally, one must evaluate
\begin{equation}
    \Cs_{n+1}^k = \Cs_{n+1}^k\left(\Qs_{n+1}^k,\ t_{n+1}\right).
\end{equation}
This vector of coefficients will not be employed for the construction of an affine source operator $\Ss_{n+1}^k$, 
but only for checking the validity of the solution obtained from the approximate problem~\eqref{eq:approxcauchy}
by comparing the coefficients vector $\Cs_{n+1}^k$ to $\Cs_{n}$, as well as comparing the coefficients $\Cs_{n+1/2}^k$
used in the middle-point affine operator for the initial coefficients $\Cs_{n}$. At the end of the iterative procedure, 
one will set $\Cs_{n+1} = \Cs_{n+1}^k$, so that this will be the new reference vector of coefficients 
for the next timestep.
The convergence criterion for stopping the iterations
is implemented by computing 
\begin{equation} \label{eq:convergencemetric}
r = \max\left(\frac{\abs{\Qs_{n+1}^{k} - \Qs_{n+1}^{k-1}}}{\abs{\Qs_{n+1}^k} + \abs{\Qs_{n+1}^{k-1}} + \epsilon_r}\right), 
\end{equation}
and checking if $r \leq r_\up{max}$, with $r_\up{max}$ and $\epsilon_r$ given tolerances, 
or if the iteration count $k$ has reached a fixed maximum value $k_\up{max}$.
Note that in principle any norm may be used to compute the error metric given in equation~\eqref{eq:convergencemetric}, as this is
just a measure of the degree to which $\Qs_{n+1}^k$ was corrected in the current iteration. Moreover, we
found convenient to limit the maximum number of iterations allowed, 
and specifically here we set $k_\up{max} = 8$, but stricter bounds can be used. 
For safety, we decide to flag the state vector $\Qs_{n+1}^{k}$ as not admissible, 
as if a floating-point exception had been triggered,
whenever the iterative procedure terminates by reaching the maximum iteration count.

After the convergence has been obtained, in order to test if the IVP \eqref{eq:exactcauchy}
is well approximated by its linearised version \eqref{eq:approxcauchy}, we compute
\begin{align}
&\delta_{n+1/2} = \max\left(\frac{\abs{\Cs_{n+1/2} - \Cs_{n}}}{\abs{\Cs_{n+1/2}} + \abs{\Cs_{n}} + \epsilon_\delta}\right), \label{eq:deltaa}\\
&\delta_{n+1} = \max\left(\frac{\abs{\Cs_{n+1} - \Cs_{n}}}{\abs{\Cs_{n+1}} + \abs{\Cs_{n}} + \epsilon_\delta}\right), \label{eq:deltab}
\end{align}
and we verify if 
$\delta = \max(\delta_{n+1/2},\ \delta_{n+1}) \leq \delta_\up{max}.$
The user should specify a tolerance $\delta_\up{max}$ as well as the floor value $\epsilon_\delta$, 
which is used in order to prevent that excessive precision requirements be imposed in those situations
when all the coefficients are so small than even large relative variations expressed by equations \eqref{eq:deltaa} and \eqref{eq:deltab}
do not affect the solution in a significant manner.
If $\delta \leq \delta_\up{max}$ we confirm the state vector at the new time level to be
$\vec{V}_{n+1} = \Qs_{n+1}^k$ and a new timestep size is computed as
\begin{equation}
\Delta t_{n+1} = \lambda\,\frac{\delta_\up{max}}{\delta + \epsilon},\quad \text{with}\quad \lambda = 0.8,\quad\epsilon = 10^{-{14}},
\end{equation}
otherwise the solution of the IVP~\eqref{eq:approxcauchy} is attempted again with a reduced timestep size, specifically one 
that is obtained by halving the timestep used in the current attempt. The same happens if
at any time the admissibility test on $\Qs_{n+1}^k$ fails.

\subsection{Analytical solution of the linearised problem}

The general solution to an initial value problem like \eqref{eq:approxcauchy} can be written as 
\begin{equation}
\vec{V}(t) = \exp\left[\Js\,(t - t_n)\right]\,\left[\vec{V}(t_n) + 
    \Js^{-1}\,\Bs - \Vs\right] - \Js^{-1}\,\Bs + \Vs.
\end{equation}
Note that, in addition to evaluating the matrix exponential $\exp\left[\Js\,(t - t_n)\right]$, one must also
compute the inverse Jacobian matrix $\Js^{-1}$. Computation of matrix exponentials can be 
carried out rather robustly in double precision arithmetic with the aid of the algorithms of 
Higham \cite{matexp2005} and Al-Mohy and Higham \cite{matexp2009, matexp2011}, while inversion of
the Jacobian matrix can be an arbitrarily ill-conditioned problem, to be carefully treated or avoided if possible.

For this reason we propose the following strategy for choosing a more suitable linearisation and 
computing analytical solutions of the linearised problem for the ODE system \eqref{eq:bn7odeu1}--\eqref{eq:bn7odea1}.
First, it is easy to see that the velocity sub-system (equations for $u_1$ and $u_2$) can be fully decoupled
from the other equations, as the partial densities $\alpha_1\,\rho_1$ and $\alpha_2\,\rho_2$ remain 
constant in the relaxation step.
\begin{figure}[!t]
\includegraphics[scale=1.05]{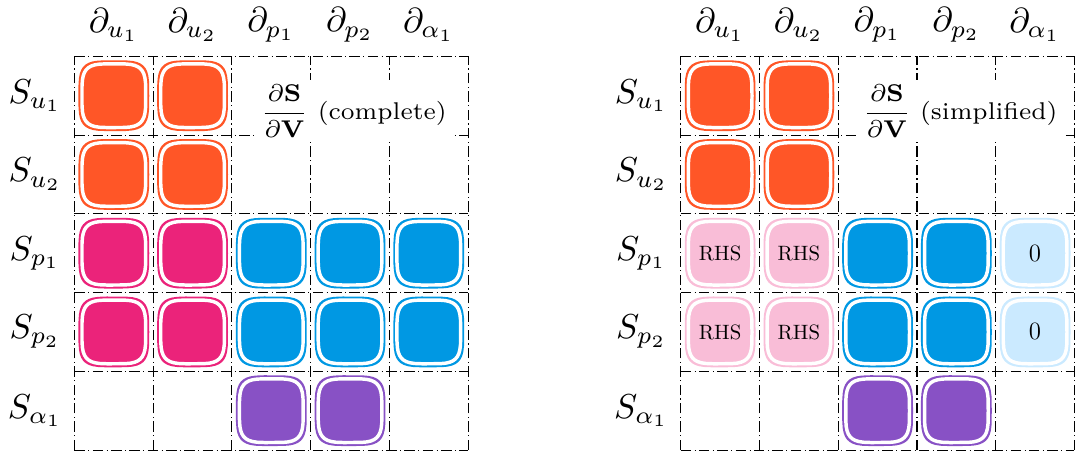}
\caption{Visual comparison between the structure of the complete Jacobian 
matrix for the ODE system \eqref{eq:bn7odeu1}--\eqref{eq:bn7odea1} and the proposed three-step simplified structure. 
The RHS label indicates dependencies that are accounted for as non-homogeneous terms in the pressure
sub-system, while the zeros mark dependencies that are suppressed entirely.}
\label{fig:Jacobian}
\end{figure}%
Then the solution of the velocity sub-system can be immediately obtained as 
\begin{align}
    & u_1(t) = \frac{\lambda}{k}\,\left\{\frac{u_1^0}{\alpha_2\,\rho_2} + \frac{u_2^0}{\alpha_1\,\rho_1} + \frac{u_1^0 - u_2^0}{\alpha_1\,\rho_1}\,\exp{\left[-k\,(t - t_n)\right]}\right\},\\
    & u_2(t) = \frac{\lambda}{k}\,\left\{\frac{u_1^0}{\alpha_2\,\rho_2} + \frac{u_2^0}{\alpha_1\,\rho_1} + \frac{u_2^0 - u_1^0}{\alpha_2\,\rho_2}\,\exp{\left[-k\,(t - t_n)\right]}\right\},
\end{align}
with $k = {1}/{\alpha_1\,\rho_1} + {1}/{\alpha_2\,\rho_2}$.
\begin{figure}[!b]
\includegraphics[scale=1.05]{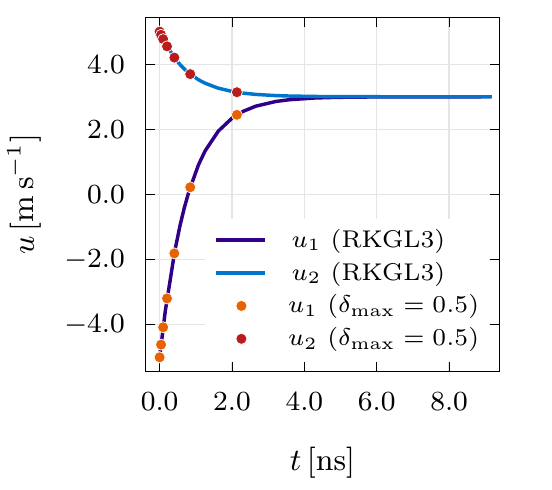}
\includegraphics[scale=1.05]{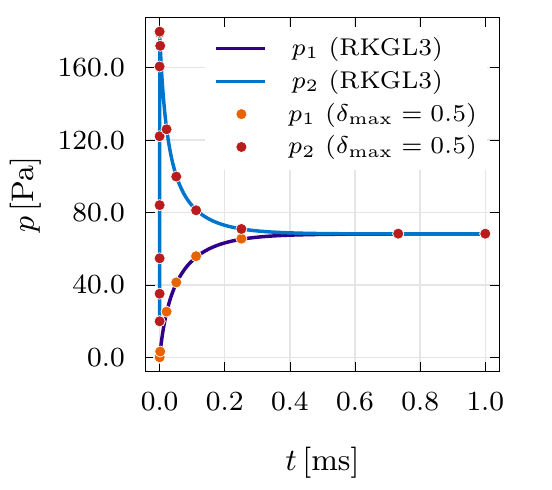}\\[3mm]
\includegraphics[scale=1.05]{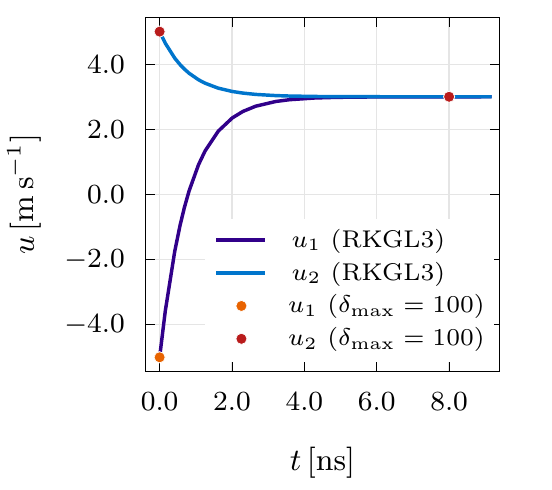}
\includegraphics[scale=1.05]{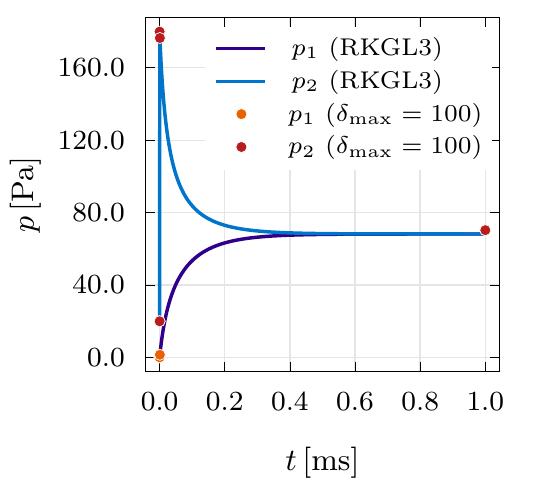}
\caption{Time evolution of velocities and pressures for test problem A1. In the top frames the linearisation tolerance
parameter is set as $\delta_\up{max} = 0.5$, employing 15 timesteps to reach the final time of $1.0\,\up{ms}$, while in the bottom frames we 
impose an extremely loose tolerance $\delta_\up{max} = 100$, still showing good agreement with the reference solution but using only 4 timesteps for
the full run.}
\label{fig:bntest1}
\end{figure}
In a second step, the pressure sub-system \eqref{eq:bn7odep1}--\eqref{eq:bn7odep2} is linearised as
\begin{align}
            & \od{p_1}{t} = k_p\,(p_2 - p_1) + k_u\,{(u_I - u_1)}\,(u_2 - u_1), \label{eq:linp1}\\[2mm]
            & \od{p_2}{t} = k_p\,(p_1 - p_2) + k_u\,{(u_I - u_2)}\,(u_1 - u_2), \label{eq:linp2}
\end{align}
where $k_p$ and $k_u$ are constant coefficients directly obtained 
from equations \eqref{eq:bn7odep1}--\eqref{eq:bn7odep2}. 
This way, at the cost of
suppressing the dependence on $\alpha_1$ in the Jacobian of the pressure sub-system, the homogeneous 
part of equations \eqref{eq:linp1}--\eqref{eq:linp2} has the same simple structure found in the velocity
sub-system, with the addition of a non-homogeneous term, which is known, 
as $u_1(t)$ and $u_2(t)$ already have been computed. The solution can again be evaluated using standard scalar exponential 
functions, which are fast and robust, compared to matrix exponentials and especially so, because one no longer
needs to perform the inversion of the Jacobian matrix of the full system.
\begin{figure}[!b]
\includegraphics[scale=1.05]{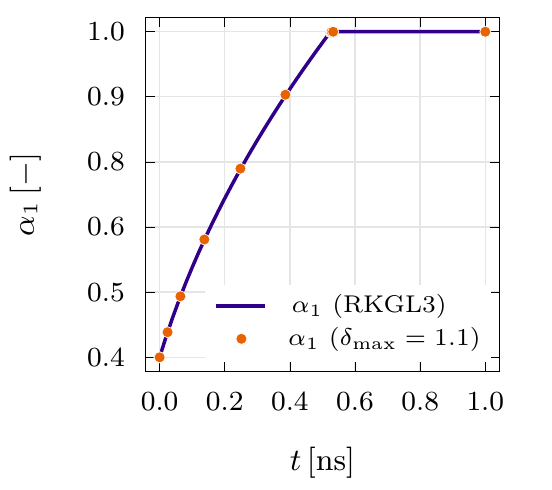}
\includegraphics[scale=1.05]{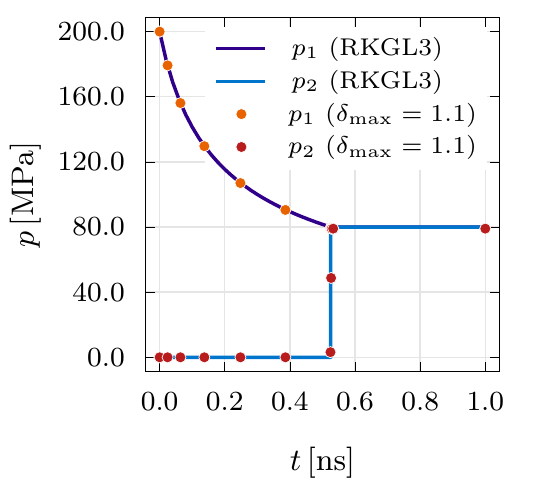}
\caption{Time evolution of volume fraction and pressure for test problem A2.
The solution is well captured in 11 timesteps, using a linearisation tolerance $\delta_\up{max} = 1.1$.}
\label{fig:bntest3}
\end{figure}
\begin{figure}[!b]
\includegraphics[scale=1.05]{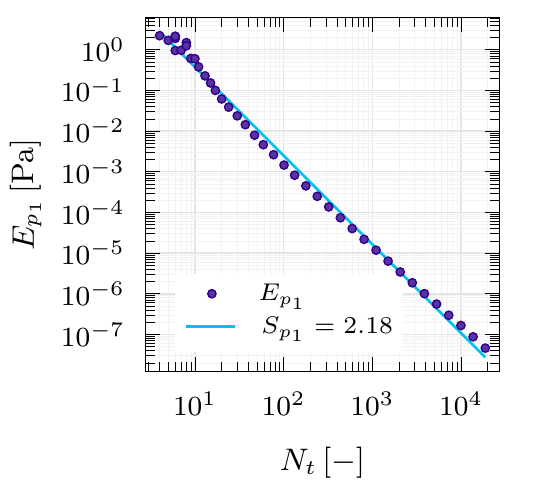}
\includegraphics[scale=1.05]{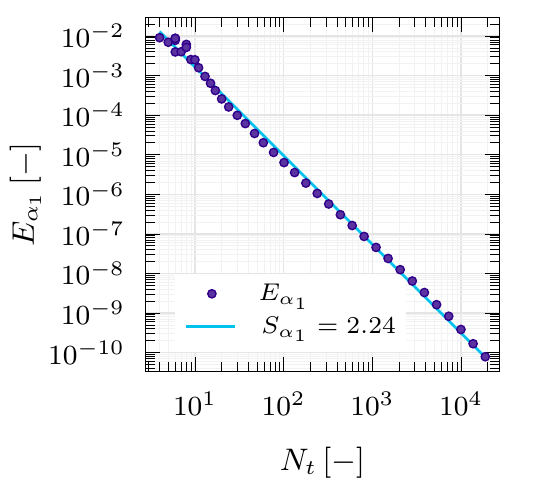}
\caption{Convergence results relative to 40 runs of test problem A1. 
On the bilogarithmic plane, the slopes of
the regression lines are $S_{p_1} = 2.18$ and $S_{\alpha_1} = 2.24$ 
for the variables $p_1$ and $\alpha_1$ respectively, indicating second order convergence.}
\label{fig:convergence}
\end{figure}
Finally, the solution to equation \eqref{eq:bn7odea1} can be integrated analytically from the expressions of
$p_1(t)$ and $p_2(t)$.
Full coupling of the system is restored in the successive iterations by 
recomputing the constant coefficients $k_p$ and $k_u$ using an updated midpoint value for $\alpha_1$.
See Figure~\ref{fig:Jacobian} for a graphical description of the proposed simplified solution structure.

\section{Test problems}
\label{sec:testproblems}
\begin{figure}[!b]
\includegraphics[scale=1.05]{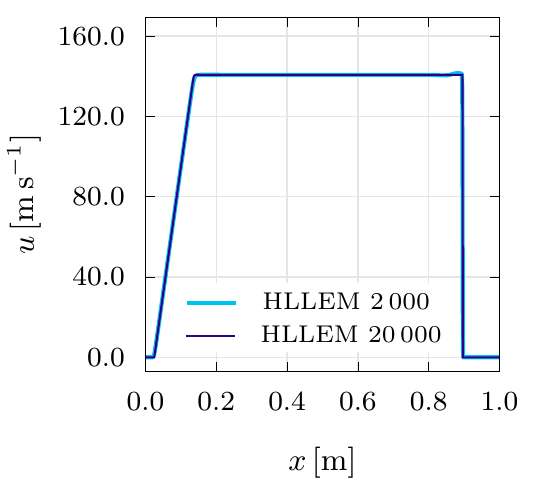}
\includegraphics[scale=1.05]{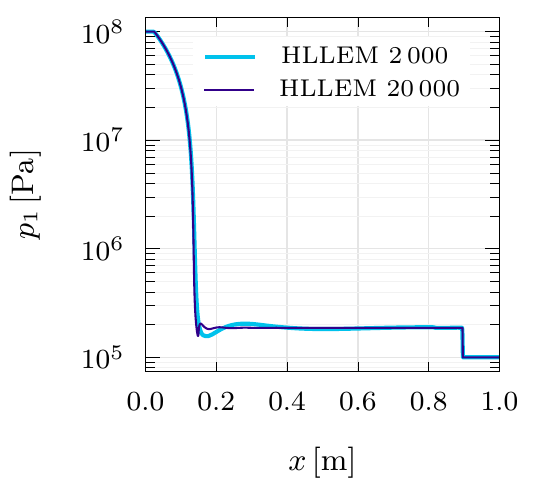}
\caption{Solution of test Problem RP1 on two uniform meshes of 2\,000 cells and 20\,000 cells respectively, 
showing convergence with respect to mesh refinement.}
\label{fig:pelanti1}
\end{figure}
\begin{figure}[!b]
\includegraphics[scale=1.05]{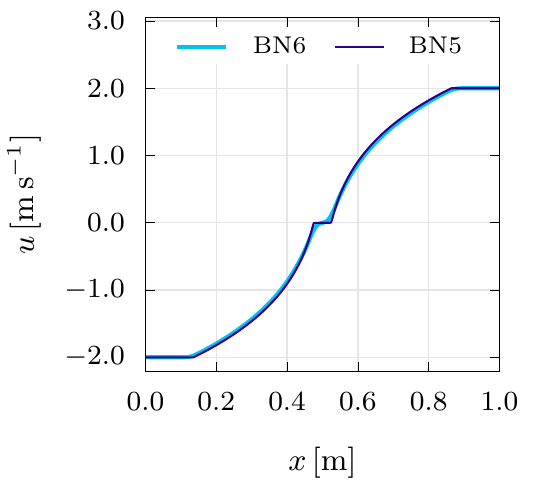}
\includegraphics[scale=1.05]{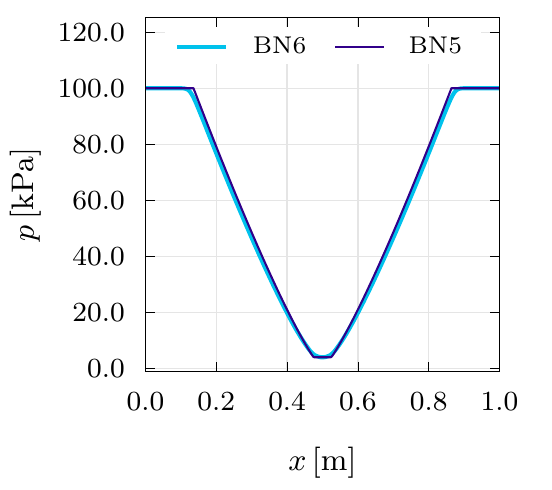}
\caption{Solution of test Problem RP2 computed from the six-equation Baer--Nunziato model (BN6) with stiff relaxation, 
compared with the five-equation Kapila model (BN5), showing convergence to the limit reduced model.}
\label{fig:pelanti2}
\end{figure}
We provide validation of the proposed method first by computing solutions 
to the ODE system \eqref{eq:bn7odeu1}--\eqref{eq:bn7odea1} and comparing the results with a reference solution
obtained from a sixth order, fully implicit, Runge--Kutta--Gauss--Legendre method \cite{butcher1964}
(labeled RKGL3) employing adaptive timestepping (test problems A1 and A2, Figures~\ref{fig:bntest1}~and~\ref{fig:bntest3}). 
Furthermore, test problem A1 is employed also for carrying out a convergence study of the scheme (Figure~\ref{fig:convergence}), 
showing that second order convergence is easily achieved.
The initial data for the ODE tests are, 
for test A1, 
\begin{equation}
u_1^0 = -5\,\up{m\,s^{-1}},\ \  u_2^0 = 5\,\up{m\,s^{-1}},\ \  p_1^0 = 0.1\,\up{Pa},\ \ p_2^0 = 20\,\up{Pa},\ \ \alpha_1^0 = 0.9,
\end{equation}
while for test A2, 
\begin{equation}
u_1^0 = 0\,\up{m\,s^{-1}},\ \  u_2^0 = 0\,\up{m\,s^{-1}},\ \  p_1^0 = 2.0\times10^8\,\up{Pa},\ \  p_2^0 = 1\,\up{Pa},\ \  \alpha_1^0 = 0.4. 
\end{equation}
The parametric data are, 
for test A1, 
\begin{equation}
\begin{aligned}
&\alpha_1\,\rho_1 = 1.0\,\up{kg\,m^{-3}},\ \ \alpha_2\,\rho_2 = 4.0\,\up{kg\,m^{-3}},\ \ \gamma_1 = 6,\ \  \gamma_2 = 1.4,\\
&\Pi_1 = 0\,\up{Pa},\ \  \Pi_2 = 0\,\up{Pa},\ \  
\lambda=10^9\,\up{kg\,m^{-1}\,s^{-1}},\ \  \nu=10\,\up{Pa^{-1}\,s^{-1}}.\\
\end{aligned}
\end{equation}
and for test A2, 
\begin{equation}
\begin{aligned}
&\alpha_1\,\rho_1 = 780.0\,\up{kg\,m^{-3}},\ \  \alpha_2\,\rho_2 = 0.22\,\up{kg\,m^{-3}},\ \  \gamma_1 = 6,\ \  \gamma_2 = 1.4,\\
&\Pi_1 = 100\,\up{Pa},\ \  \Pi_2 = 0\,\up{Pa},\ \ 
\lambda=10^9\,\up{kg\,m^{-1}\,s^{-1}},\ \  \nu=10\,\up{Pa^{-1}\,s^{-1}}.
\end{aligned}
\end{equation}
Then, we show an application
of the method
in the solution of the mixture-energy-consistent formulation of 
the six-equation reduced Baer--Nunziato model forwarded in \cite{pelanti}.
For these
simulations the interface pressure is computed as 
\begin{equation}
    p_I = \frac{Z_2\,p_1 + Z_1\,p_2}{Z_1 + Z_2},\quad \text{with}\quad Z_1 = \rho_1\,a_1\ \text{and}\ Z_2 = \rho_2\,a_2.
\end{equation}

\begin{figure}[!p]
\includegraphics[scale=1.00]{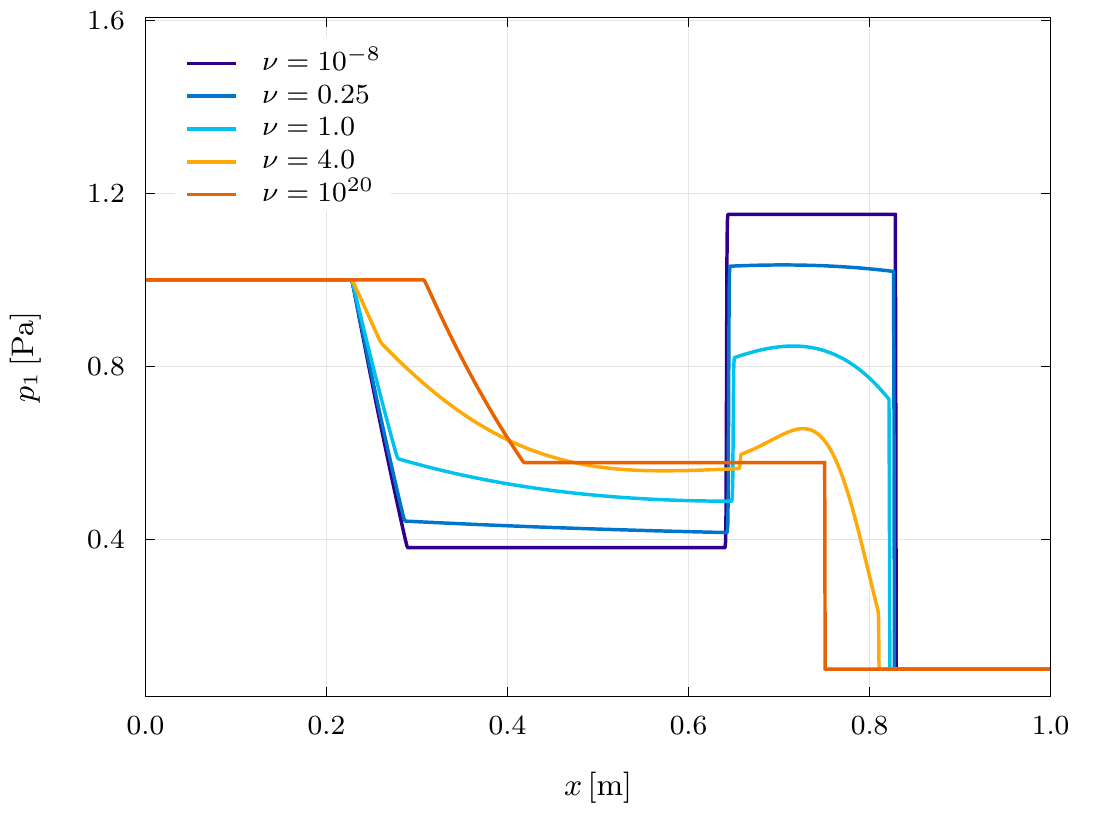}\\[4mm]
\includegraphics[scale=1.00]{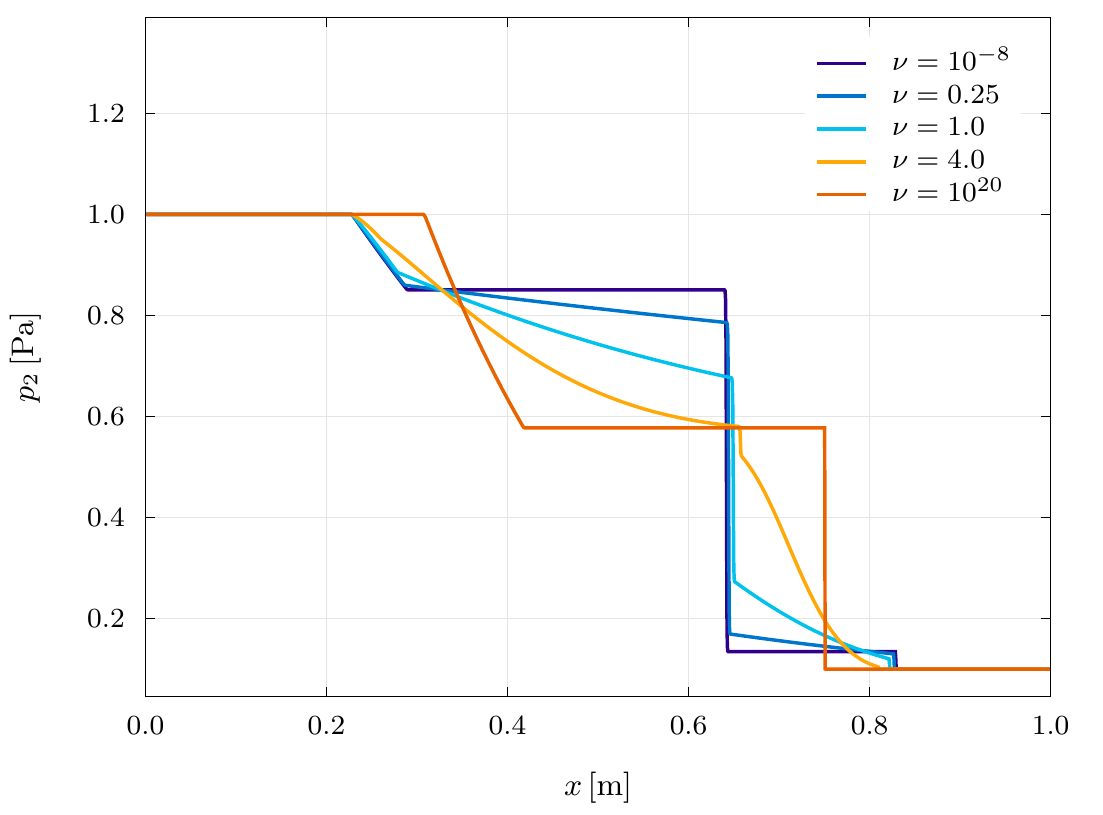}
\caption{Behaviour of the pressure variables in RP3 with several values of $\nu$. It is clear that, in the 
stiff regime ($\nu = 10^{20}\,\up{Pa^{-1}\,s^{-1}}$), $p_1$ and $p_2$ converge to the same value, while they evolve in a 
completely distinct fashion if relaxation is set to act on longer timescales.}
\label{fig:finrate12}
\end{figure}

The first two shock-tube problems (from \cite{cavitation_shocktube, pelanti}), 
show that the method is able to deal with very stiff ($\nu = 10^{20}\,\up{Pa^{-1}\,s^{-1}}$) sources, and 
in particular in Figure~\ref{fig:pelanti1} (RP1, a liquid-vapour dodecane shock tube 
featuring a strong right-moving shockwave) 
we show mesh convergence of the solution by comparing two runs, 
both employing the HLLEM Riemann solver proposed in \cite{hllem}, on two different meshes consisting of
2\,000 uniform control volumes and 20\,000 control volumes respectively, 
with a computational domain delimited by $x \in [0\,\up{m}, 1\,\up{m}]$.
In Figure~\ref{fig:pelanti2} (RP2, two diverging rarefaction waves in liquid water) we then show that, 
with very stiff relaxation ($\nu = 10^{20}\,\up{Pa^{-1}\,s^{-1}}$), 
the solution matches the one computed by solving directly the five-equation 
instantaneous equilibrium model \cite{kapila}, 
again using a mesh consisting of 2\,000 uniform cells for the six-equation model and
a mesh of 20\,000 uniform cells for the reference solution, and in particular, rarefaction waves
propagate at the same speed for both models.
All tests are run using a second order path-conservative MUSCL-Hancock scheme with $k_\up{CFL}=0.95$.
The first Riemann Problem (RP1) is set up with uniform liquid and vapour 
densities $\rho_1^\up{L} = \rho_1^\up{R} = 500\,\up{kg\,m^{-1}\,s^{-1}}$ and $\rho_2^\up{L} = \rho_2^\up{L} = 2.0\,\up{kg\,m^{-3}}$, 
uniform velocity $u^\up{L} =  u^\up{R} = 0\,\up{m\,s^{-1}}$, a jump in pressure given by 
$p_1^\up{L} = p_2^\up{L} = 100\,\up{MPa}$, $p_1^\up{R} = p_2^\up{R} = 100\,\up{kPa}$, 
almost pure liquid on the left side of the initial discontinuity ($\alpha_1^\up{L} = 1-10^{-8}$), 
and almost pure vapor on the right side ($\alpha_1^\up{R} = 10^{-8}$).
The discontinuity is initially found at $x = 0.75\,\up{m}$,   
and the end time is $t_\up{end} = 473\,\up{\mu s}$. The parameters of the stiffened gas EOS
are  $\gamma_1 = 2.35$, $\gamma_2 = 1.025$, $\Pi_1 = 400\,\up{MPa}$, $\Pi_2 = 0$.

The second Riemann Problem (RP2) is initialised with constant liquid and vapour densities 
$\rho_1^\up{L} = \rho_1^\up{R} = 1150\,\up{kg\,m^{-1}\,s^{-1}}$,  $\rho_2^\up{L} = \rho_2^\up{L} = 0.63\,\up{kg\,m^{-3}}$, 
constant pressure $p_1^\up{L} = p_2^\up{L} = p_1^\up{R} = p_2^\up{R} = 100\,\up{kPa}$, constant liquid volume 
fraction $\alpha_1^\up{L} = \alpha_1^\up{R} = 0.99$, and a jump in velocity (initially located at $x = 0.5\,\up{m}$) 
such that $u^\up{L} = -2.0\,\up{m\,s^{-1}}$ and $u^\up{R} = 2.0\,\up{m\,s^{-1}}$.
The final time is $t_\up{end} = 3.2\,\up{ms}$ and for this test the parameters of the equation of state 
$\gamma_1 = 2.35$, $\gamma_2 = 1.43$, $\Pi_1 = 1\,\up{GPa}$, $\Pi_2 = 0$.

Finally, in Figure~\ref{fig:finrate12} we show
the behaviour of the solution of a third Riemann problem (RP3) 
with several different values of the pressure relaxation parameter $\nu$ 
(ranging from $10^{-8}\,\up{Pa^{-1}\,s^{-1}}$ to $10^{20}\,\up{Pa^{-1}\,s^{-1}}$), highlighting 
the vast range of solution structures that can be obtained not only with stiff 
relaxation (the pressure profiles $p_1$ and $p_2$ coincide) or in total absence of it, but also 
with finite values of the relaxation time scale.
For RP3, the initial data on the left are
\begin{equation}
\begin{aligned}
    & \rho_1^\up{L} = 1.0\,\up{kg\,m^{-1}\,s^{-1}},\quad  \rho_2^\up{L} = 0.2\,\up{kg\,m^{-3}},\quad  u^\up{L} = 0.0\,\up{m\,s^{-1}},\\
    & p_1^\up{L} = 1.0\,\up{Pa},\quad  p_2^\up{L} = 1.0\,\up{Pa},\quad \alpha_1^\up{L} = 0.55,
\end{aligned}
\end{equation}
while on the right one has 
\begin{equation}
\begin{aligned}
& \rho_1^\up{R} = 0.125\,\up{kg\,m^{-3}},\quad  \rho_2^\up{R} = 2.0\,\up{kg\,m^{-3}},\quad u^\up{R} = 0.0\,\up{m\,s^{-1}},\\ 
& p_1^\up{R} = 0.1\,\up{Pa},\quad  p_2^\up{R} = 0.1\,\up{Pa},\quad \alpha_1^\up{R} = 0.45. 
\end{aligned}
\end{equation}

The initial jump is located at $x = 0.6\,\up{m}$, the domain is $x \in [0\,\up{m},\ 1\,\up{m}]$ and the final time is $t_\up{end} = 0.15\,\up{s}$.
The parameters of the stiffened gas EOS are $\gamma_1 = 2.0$, $\gamma_2 = 1.4$, $\Pi_1 = 2.0\,\up{Pa}$, $\Pi_2 = 0.0\,\up{Pa}$.

\section{Conclusions}
We presented a technique for integrating ordinary differential equations associated with
stiff relaxation sources and promising results have been shown for a set
of test problems. The method can efficiently resolve very abrupt variations in the solution and adapt
to multiple timescales. A key feature of the algorithm is that it can avoid delicate linear algebra 
operations entirely, thus improving the robustness of the scheme.
Future applications will include liquid-gas and liquid-solid phase transition, strain relaxation
for nonlinear elasticity \cite{gpr} and the computation of material failure in elasto-plastic and brittle solids.

\acknowledgement{The authors of this work were supported by the German Research Foundation (DFG) 
through the project GRK 2160/1 ``Droplet Interaction Technologies''.}

This is a pre-print of the following work: 
{G. Lamanna, S. Tonini, G.E. Cossali and B. Weigand (Eds.)}, 
\textit{``Dropet Interaction and Spray Processes''}, 2020, 
Springer, Heidelberg, Berlin. 
Reproduced with permission of Springer Nature Switzerland AG.

\texttt{DOI: 10.1007/978-3-030-33338-6}.

\bibliographystyle{plain}
\bibliography{./references}

\end{document}